\documentclass[preprint]{elsarticle}

\usepackage{lineno,hyperref}

\usepackage{graphicx}
\usepackage{amsmath,amssymb}
\usepackage{verbatim}
\usepackage{mathptmx}      
\usepackage{latexsym}
\usepackage{marvosym}
\modulolinenumbers[5]

\journal{Journal of \LaTeX\ Templates}
\newproof{proof}{Proof}
\newtheorem{thm}{Theorem}
\newtheorem{lem}{Lemma}
\newdefinition{rmk}{Remark}
\newtheorem{prop}{Proposition}
\newproof{pot}{Proof of Theorem \ref{thm2}}

\newcommand{\grad}{\mathrm{grad\ }}
\newcommand{\dist}{\mathrm{dist}}
\newcommand{\supp}{\mathrm{supp}}

\begin{document}

\begin{frontmatter}

\title{ H{\"o}lder  classes  in $L^p$ norm on a chord arc curve in $\mathbb R^3$}

\author[rvt]{Tatyana~A.~Alexeeva\corref{cor1}}
\ead{tatyanalexeeva@gmail.com}

\author[rvt,focal]{Nikolay~A.~Shirokov}
\ead{nikolai.shirokov@gmail.com}

\cortext[cor1]{Corresponding author}
\address[rvt]{Department of Mathematics,  St.~Petersburg School of Physics, \\ Mathematics, and Computer Science, HSE University,

3A Kantemirovskaya Ul., St.~Petersburg, 194100, Russia}


\address[focal]{Department of Mathematical Analysis, Faculty of Mathematics and Mechanics,

St.~Petersburg State University, 28 Universitetsky prospekt, Peterhof, St.~Petersburg, 198504, Russia}




\fntext[myfootnote]{The second author was supported by the RFBR grant 20-01-00209}

\begin{abstract}
We define H{\"o}lder  classes  in the $L^p$ norm   on a chord-arc curve  in $\mathbb{R}^3$ and prove direct and inverse approximation theorems for   functions from these classes by functions harmonic in a neighborhood of the curve. The approximation is estimated in the $L^p$ norm, and the smaller the neighborhood, the more accurate the approximation.
\end{abstract}

\begin{keyword}
Constructive description\sep H{\"o}lder classes \sep Harmonic functions \sep Chord-arc curves
\MSC[2010] 41A30\sep  41A27
\end{keyword}

\end{frontmatter}

\linenumbers

\section{Introduction}
The problem of describing H{\"o}lder classes in the $L^p$ norm was first considered in the case of periodic functions in terms of the approximation rate by trigonometric polynomials \cite [Ch.6]{Timan-Book-1963}. A description of the approximation of smooth functions by algebraic polynomials in the $L^p$ norm was obtained by V.P.Motornyi \cite {Motornyi-1971} in 1971. In \cite {Potapov-1977} M.K.Potapov gave a constructive description of the approximation of new classes of functions on an interval by algebraic polynomials in terms of the approximation rate in the $L^p$ norm. These classes were defined by weight conditions in the same $L^p$ norm. The construction of approximating polynomials was nontrivial. P.Nevai and Yuan Xu \cite {NevaiYu-1994} used some other polynomials in a similar case of the $L^p$ approximation.
E.M.Dynkin \cite {Dynkin-1983} extended the description in terms of approximation by algebraic polynomials
in the $L^p$ norm to the Sobolev and Besov classes  on an interval. He also obtained a constructive description of the same classes of functions analytic in Jordan domains. The boundaries of these domains were assumed to have the property of commensurability of arcs and chords;  nowadays such boundaries are called chord-arc curves.

In E.M. Dynkin's description of the approximation in a domain, as well as in the statements of the results concerning a constructive description of classes of functions on an interval, the $L^p$ norm approximation scales connected with a conformal map from the exterior of a domain or an interval to the exterior of the unit disc were applied.  However, if we want to get a constructive description of functional classes on a curve in $\mathbb R^3$, we cannot use  constructions similar to those used in the case of the complex plane. Besides, the problem of a constructive description of functional classes on a nonclosed curve in terms of uniform approximation by polynomials  has features different from those encountered in the description of functions analytic in Jordan domains. In the case of the approximation in the uniform norm, V.V.Andrievskii \cite {Andrievskii95} used estimates for polynomial approximation rate and  for the derivatives of approximating polynomials and obtained a constructive description of the classes of smooth functions defined on a nonclosed chord-arc curve in the plane.

In the paper \cite {AlexeevaShir-2020}, the authors used an approach based on the estimation of the approximation rate and the gradient of the approximating function, which allowed them to get a constructive description of H\"{o}lder type classes of functions defined on a nonclosed chord-arc curve in $\mathbb R^3$. Approximations were obtained by   functions harmonic in a neighborhood of the curve. The neighborhoods were compressed to obtain a better approximation. For plane domains, a constructive description of functional classes by means of the approximation by harmonic polynomials was obtained by V.V.Andrievskii \cite{Andrievskii88}.

In the present paper, we consider the approximation rate in the $L^p$ norm for functions that are
defined on a chord-arc curve in $\mathbb R^3$ and belong to the classes that can be named H\"{o}lder in the $L^p$ norm. As approximating functions we take the same  harmonic functions
defined in compressing neighborhoods of the curve that were used for uniform approximation in \cite {AlexeevaShir-2020}. The theorem on a possible approximation rate is proved for a narrower class of functions than the theorem on the smoothness of a function approximated with the above-mentioned rate.

The paper is organized as follows.  Sec.\ref{sec2} contains main definitions and statement of   main results. In Sec.\ref{sec3}  we prove Theorem\ref{thm2}  saying that a function approximated with a certain rate is smooth in the $L^p$ norm. In Sec.\ref{sec4} we construct a special continuation of a function defined on a curve to the entire $\mathbb R^3$. In Sec.\ref{sec5} we construct the approximating functions. In Sec\ref{sec6} we prove Theorem\ref{thm1} on a possible approximation rate. Sec.\ref{sec7}  is a conclusion section.

\section{Definitions and statement of   main results}\label{sec2}
Let $L$ be a nonclosed curve in $\mathbb R^3$ with the endpoints $A$  and $B$.
We say that $L$ is $b$-chord-arc, where $b\geq 1$, if for all $M_1, M_2\in L$,
$M_1\neq M_2$, the arc $\gamma(M_1, M_2)\subset L$  with the endpoints $M_1$ and $M_2$ satisfies the inequality
$|\gamma(M_1,M_2)|\leq b||M_1M_2||$, where $|\gamma(M_1,M_2)|$ means the length of $\gamma(M_1,M_2)$.
Denote $B_r(M)=\{N\in\mathbb R^3:||MN||<r\}$ and $\bar B_r(M)=\{N\in\mathbb R^3:||MN||\leq r\}$. For a function $f$ defined on $L$ we put
$$\Delta^{\ast}f(M,r)=\underset{N\in \bar B_r(M)\cap L}\sup|f(N)-f(M)|, \ M\in L.$$
For $0<\alpha<1$ we denote by $\Lambda^{\alpha}_p(L)$ the space of all functions defined on $L$ and satisfying the condition
\begin{equation}\label{1}\underset{0<r<|L|}\sup\left(\int\limits_{L}
\left(\frac{\Delta^{\ast}f(M,r)}
{r^{\alpha}}\right)^p\,dm_1(M)\right)^{\frac1p}<\infty,
\end{equation}
where $m_1(M)$ is the curve length on $L$.
Let $\Lambda^{\alpha_0}_p(L)$ be the subspace of $\Lambda^{\alpha}_p(L)$ consisting of the functions for which condition \eqref{1} is fulfilled and for some $\epsilon=\epsilon(f)>0$ and $c=c(f)>0$ the following inequality is valid:
\begin{equation}\label{2}
\Delta^{\ast}f(M,r)\leq c\left(\frac rR\right)^{\epsilon}\Delta^{\ast}f(N,R),
\end{equation}
where $0<r\leq R$,\ $||MN||\leq R$,\,  $M,  N\in L$.
Let $\Lambda^{\alpha}$ be the space of functions $f$ on $L$ for which
$$\Delta^{\ast}f(M,r)\leq cr^{\alpha}, \ M\in L,\ \text  {and}\;\  c=c(f).$$
Let $\Omega_{\delta}=\underset{M\in L}\bigcup B_{\delta}(M)$ and let $H(\Omega_{\delta})$ be the set of functions harmonic in $\Omega_{\delta}$. For
$v\in H(\Omega_{\delta})$ we put
$$\grad_{\delta}^{\ast}v(M)=\underset{N\leq \bar B_{\frac{\delta}{2}}(M)}\max|\grad v(N)|,$$
and for a function $F$ defined on $L$ we put
$$\max\nolimits_{\delta}F(M)=\underset{N\in\bar B_{\delta}(M)\cap L}\sup|F(N)|.$$
The following statements are valid.
\begin{thm}\label{thm1}
Let $f\in\Lambda_p^{\alpha_0}(L)$, $0<\alpha<1$, $p>1/\alpha$. Then $f\in\Lambda^{\alpha-\frac1p}(L)$ and there is a constant $c=c(f)>0$ such that for each $0<\delta<|L|$ there exists a function $v_{\delta}\in H(\Omega_{\delta})$ for which the following inequalities are valid:
\begin{equation}\label{3}
\left(\int\limits_L\left(\frac{\max_{\delta}(f(M)-v_{\delta}(M))}
{\delta^{\alpha}}\right)^p\,dm_1(M)\right)^{\frac1p}\leq c
\end{equation}
and
\begin{equation}\label{4}
\left(\int\left(\delta^{1-\alpha}\grad^{\ast}_{\delta}\ v_{\delta}(M)\right)^p\,dm_1(M)\right)^{\frac1p}\leq c.
\end{equation}
\end{thm}
If a function $f$ can be approximated so that conditions \eqref{3} and \eqref{4} are satisfied, then $f\in \Lambda^{\alpha}_p(L)$.

\begin{thm}\label{thm2}
Let $f\in C(L)$ and for $0<\delta\leq 2|L|$ there exist functions $v_{\delta}\in H(\Omega_{\delta})$ satisfying conditions \eqref{3} and \eqref{4}, where $0<\alpha<1$ and
$p>1/\alpha$. Then $f\in\Lambda_{p}^{\alpha}(L)$.

\end{thm}
\section{Proof of Theorem \ref{2}}\label{sec3}
Let $N(M)$ be a function measurable with respect to the $m_1$-Lebesgue measure on $L$ such that $||MN(M)||\leq r$, $0<r\leq |L|$. We choose $c_1>2$   such that $[M,N(M)]\subset \Omega_{\frac{c_1r}  {2}}(L)$ for $M\in L$ and
$0<r\leq |L|$.
Since $L$ is a $b$-chord-arc curve, such a  choice is possible and is described in detail in \cite {AlexeevaShir-2020}.
Denote by $\nu(M)$ a unit vector parallel to the vector $\overline{MN}(M)$ if $n(M)\neq M$; the vector $\nu(M)$ is undefined if $N(M)=M$.
We put $\delta=c_1r$ and choose a function $v_{\delta}\in H(\Omega_{\delta}(L))$ satisfying conditions \eqref{3} and \eqref{4}. Then we have
\begin{equation}\label{5}
\begin{aligned}f(N(M))-f(M)=(f(N(M))-v_{\delta}(N(M)))-(f(M)-v_{\delta}(M))+
v_{\delta}(N(M))-v_{\delta}(M)\\=
((f(N(M))-v_{\delta}(M))-(f(M)-v_{\delta}(M))+\int\limits_0^{||MN(M)||}
v^{\prime}_{\delta\nu(M)}(M+t\nu(M))\,dt,
\end{aligned}\end{equation}
where $v^{\prime}_{\delta\nu(M)}(M+t\nu(M))$ is the directional derivative along the vector $\nu(M)$.
Now \eqref{5} implies the estimate
\begin{equation}\label{6}\begin{aligned}
|f(N(M))-f(M)|\leq\max\nolimits_{\delta}(f(M)-v_{\delta}(M))+
\max\nolimits_{\delta}(f(M)-v_{\delta}(M))\\+
||MN(M)|| \grad_{\delta}^{\ast}v(M)\leq 2\max\nolimits_{\delta}(f(M)-v_{\delta}(M))+
\delta \grad_{\delta}^{\ast}v(M).
\end{aligned}
\end{equation}
From \eqref{6}
we obtain that
\begin{equation}\label{7}
\begin{aligned}
\left(\int\limits_L\left|\frac{f(N(M))-f(M)}{r^{\alpha}}\right|^p dm_1(M)\right)^{\frac1p}
\leq  c \left(\int\limits_L
\left(\frac{\max_{\delta}(f(M)-v_{\delta}(M))}{\delta^{\alpha}}\right)^p
dm_1(M)\right)^{\frac1p}\\
+\left(\int\limits_L(\delta^{1-\alpha}\grad_{\delta}^{\ast}v(M))^p
dm_1(M))^{\frac1p}\right)\leq c^{\prime},
\end{aligned}
\end{equation}
where the constant $c^{\prime}$ does not depend on the choice of $N(M)$
and $r>0$. Due to an arbitrary choice of the function $N(M)$ satisfying the condition
$||MN(M)||\leq r$, estimate \eqref{7} implies \eqref{1}. Theorem~2 is proved.
\section{Pseudoharmonic continuation of the function f}\label{sec4}
\begin{prop}\label{prop1}
Let $f$ be a function satisfying \eqref{1}. Then $f\in\Lambda^{\alpha-\frac1p}(L)$.
\end{prop}
\begin{proof}
Let $M, M_1, M_2\in L$ be such that the point $M$ lies on the curve $\gamma(M_1,M_2)$,
$|\gamma(M_1,M)|=|\gamma(M,M_2)|=r$, $0<r<L$, and $M_1=A$ if
$|\gamma(A,M)|<r$ and $M_2=B$ if $|\gamma(M,B)|<r$.
Then $||M_1M||\leq r$ and $||MM_2||\leq r$.
Notice that for $N\in \gamma(M_1,M_2)$ we have the inequality
$\Delta^{\ast}f(M,r)\leq2\Delta^{\ast} f(N,2r)$. Therefore,
\begin{equation}\label{8}
\begin{aligned}
\int\limits_{\gamma(M_1,M_2)}
\left(\frac{\Delta^{\ast}f(M,r)}{(2r)^{\alpha}}\right)^p d m_1(N)&\leq
\int\limits_{\gamma(M_1,M_2)}
\left(\frac{2\Delta^{\ast}f(N,2r)}{(2r)^{\alpha}}\right)^p d m_1(N)\\
&\leq \int\limits_{L}
\left(\frac{2\Delta^{\ast}f(N,2r)}{(2r)^{\alpha}}\right)^p d m_1(N)\leq c,
\end{aligned}
\end{equation}
where the constant $c$ does not depend on $r$. From \eqref{8} we obtain
$$r^{1-\alpha p}(\Delta_{\ast}f(M,r))^p\leq c^{\prime},\
  \Delta^{\ast}f(M,r)\leq c^{\prime\frac1p}r^{\alpha-\frac1p},$$
which proves the inclusion $f\in \Lambda^{\alpha-\frac1p}(L)$.
\end{proof}
Let us construct a continuation of $f$ to the entire space $\mathbb R^3$ in the same way as in \cite{AlexeevaShir-2020}.
Let $n\geq 1$. We subdivide the curve $L$ into $2^n$ parts of equal length by the points
$M_{kn}$, where $M_{0n}=A$, $M_{2^n,n}=B$, $\gamma(M_{0n},M_{kn})\subset
\gamma(M_{0n},M_{k+1,n})$, $1\leq k
\leq 2^n-1$.
Let
$$\begin{aligned}\Lambda_n&=|\gamma(M_{kn},\ M_{k+1,n})|=2^{-n}|L|, \
\Omega_n=\Omega^{\ast}_n\setminus \Omega^{\ast}_{n+1},\
\omega_{0n}=B_{2\Lambda_n}(M_{0n})\cap\Omega_n,\\
\omega_{kn}&=(B_{2\Lambda_n}(M_{kn})\cap \Omega_n)\setminus\bigcup\limits_{\nu=0}^{k-1}B_{2\Lambda_n}(M_{\nu n}), 1\leq k \leq2^n.
\end{aligned}
$$
The set $\{\omega_{kn}\}$ may be empty for some $k$ and $n$.
For $M\not\in L$  we put
\begin{equation}\label{9}
f_1=\begin{cases}
f(M_{kn}), \ M\in \omega_{kn},\\
0, \ M\in\mathbb R^3\setminus\bigcup\limits_{n=0}^{\infty}\Omega^{\ast}_n.
\end{cases}
\end{equation}
Denote $d(M)=\dist(M,L)$ for $M\in\mathbb R^3\setminus L$.
We will use the following statement proved in \cite{AlexeevaShir-2020} (another construction is given in \cite [Ch.6]{Stein-B-1970}
\begin{prop}\label{prop2}
There exists a function $d_0(M)\in C^2(\mathbb R^3\setminus L)$ with the following properties:
\begin{equation}\label{10}
c_2d(M)\leq d_0(M)\leq\frac1{16}d(M),
\end{equation}
\begin{equation}\label{11}
||\grad d_0(M)||\leq c_3, \  M\in\mathbb R^3\setminus L,
\end{equation}
\begin{equation}\label{12}
||\grad^2d_0(M)||\leq c_4 d^{-1}(M),
\end{equation}
where the constants $c_2,c_3,c_4$ in \eqref{10} -- \eqref{12} depend on $b$.
\end{prop}
Now we put
\begin{equation}\label{13}
f_2(M)=\frac1{|B_{d_0(M)}(M)|}\int\limits_{B_{d_0(M)}(M)}f_1(K)dm_3(K),
\end{equation}
\begin{equation}\label{14}
f_0(M)=\frac1{|B_{d_0(M)}(M)|}\int\limits_{B_{d_0}(M)}f_2(K)dm_3(K).
\end{equation}
In formulas \eqref{13} and \eqref{14}, $|B_{d_0(M)}(M)|$ means the three-dimensional measure of the ball $B_{d_0(M)}(M)$. The construction of the functions $f$, $f_1$,   $f_2$ and a reasoning similar to that
in \cite {AlexeevaShir-2020}, which uses properties \eqref{10} -- \eqref{12} of the function $d_0(M)$,
imply the following result.
\begin{lem}\label{lem1}
There exists an absolute constant $c_5$ such that for $M\in B_{2\Lambda_n}(M_0)$,  $M_0\in L$,
the following estimates are valid:
\begin{equation}\label{15}
||\grad f_0(M)||\leq c\Lambda_n^{-1}\Delta^{\ast}f(M_0,c_5 \dist(M,L)),
\end{equation}
\begin{equation}\label{16}
||\grad^2 f_0(M)||\leq c\Lambda_n^{-2}\Delta^{\ast}f(M_0,c_5 \dist(M,L)),
\end{equation}
where the constant $c$ in \eqref{15}, \eqref{16} depends on $f$ and the constant
$c_5$ depends on $b$.
\end{lem}
Notice that under conditions of Lemma~1, estimate \eqref{15} and Proposition~1
yield the estimate
$$||\grad f_0(M)||\leq c\dist^{\alpha-\frac1p-1}(M,L).$$
Let $n\geq 2$,\ $m\geq n$,\ $0\leq k_0\leq 2^n$,\ $k_1=\max(k_0-2,0)$,\
$k_2=\min(k_0+2,2^n)$. We put
\begin{equation}\label{17}
\Omega_m^{\ast}(2^{m-n}k_1,2^{m-n}k_2)=\bigcup\limits_{2^{m-n}k_1}^{2^{m-n}k_2}
\bar B_{\Lambda_m}(M_{jm}).\end{equation}
For $M\in\Omega_m^{\ast}(2^{m-n}k_1,2^{m-n}k_2)$ we define
\begin{equation}\label{18}
\tilde d(M)=\underset{N\in\gamma(M_{2^{m-n}k_1,m},M_{2^{m-n}k_2,m})}\min||MN||.
\end{equation}
Let $N(M)\in \gamma(M_{2^{m-n}k_1,m},M_{2^{m-n}k_2,m})$ be any point
for which the equality in \eqref{18} is attained.
\begin{lem}\label{lem2}
Let $f\in\Lambda_p^{\alpha_0}(L)$, let $c$ and $\epsilon$  be the constants from \eqref{2}, and let $c_6\geq 4$.
Then there exists a constant $c_7=c_7(c,\epsilon, b)>0$ such that the following estimate is valid:
$$\underset{\Omega^{\ast}_m(2^{m-n}k_1,2^{m-n}k_2)}
\int\frac{\Delta^{\ast}f(N(M), c_6\tilde d(M))}{\tilde d^2(M)}dm_3(M)\leq
c_7\Lambda_n^{1-\epsilon}\Lambda_m^{\epsilon}\Delta^{\ast}f(M_{k_0n},(c_6+3b)\Lambda_n).$$
\end{lem}
\begin{proof}
We put $\Omega^0=\Omega^{\ast}(2^{m-n}k_1,2^{m-n}k_2)$,
$k_1^0=2^{m-n}k_1$, $k_2^0=2^{m-n}k_2$. Choose a $\nu$ such that
$\lambda\overset{\mathrm{def}}=b2^{-\nu}\leq 1/8$, and let
$\Omega^1=\Omega^0\cap\Omega^{\ast}_{m+\nu}$.
For a point $M\in\Omega^1$ let $k$ be such that $k_1^0\leq k\leq k_2^0$ and
$M\in \bar B_{\Lambda_m}(M_{km})$, and let $P_0(M)\in L$ be a point
such that $||MP_0(M)||=\dist(M,L)$. Then $||MP_0(M)||\leq2\Lambda_m$, $M_{km}P_0(M)||\leq4\Lambda_m$. From the properties of the curve $L$
we obtain $|\gamma(M_{km}P_0(M))|\leq4b\Lambda_m$. Let $a=[4b]+1$, $k_1^1=\max(k_1^0-a,0)$, $k_2^1=\min(k_2^0+a,2^m)$.  We put $\Omega^2=\Omega^1\cap\Omega_{m+2\nu}^{\ast}$. For
$M\in\Omega^2$ we denote by $M_{j,m+\nu} $ a point
for which $M\in\bar B_{\Lambda_{m+\nu}}(M_{j,m+\nu})$ and by
$P_1(M)\in L$ a point for which $||MP_1(M)||=\dist(M,L)$.
We take  into account that $2^{\nu}k_1^1\leq j\leq 2^{\nu}k_2^1$ and that $||MP_1(M)||\leq2\Lambda_{m+\nu}||$,
$||M_{j,m+\nu}P_1(M)||\leq 4\Lambda_{m+\nu}$,
$|\gamma(M_{j,m+\nu},P_1(M)|\leq4b\Lambda_{m,\nu}$. We put
$\Omega^3=\Omega^2\cap\Omega^{\ast}_{m+3\nu}$, etc.; if $k_1^q$ and $k_2^q$ have already been chosen and $\Omega^q\subset\Omega^{\ast}_{m+q\nu}$, then we put
$k_1^{q+1}=\max(0,2^{\nu}k_1^q-a)$, $k_2^{q+1}=\min(2^{m+q\nu},2^{\nu}k_2^q+a)$,
$\Omega^{q+1}=\Omega^q\cap\Omega^{\ast}_{m+(q+1)\nu}$.
If $M\in\Omega^q$, then $M\in\bar B_{\Lambda_{m+q\nu}}(M^q_{k,m+q\nu})$, where
$2^{\nu}k_1^q\leq k\leq 2^{\nu}k_2^q$. If $M\in\Omega^{q+1}$, then
let $P_q(M)\in L$ be such that $||MP_q(M)||=\dist(M,L)$.
Then
$||M^q_{k,m+q\nu}P_q(M)||\leq 4\Lambda_{m+q\nu}$ and
$|\gamma(M^q_{k,m+q\nu},P_q(M)|\leq4b\Lambda_{m+q\nu}$.
From the obtained estimates we find that for $q\geq1$ we have
\begin{equation}\label{19}
|\gamma(M^{q+1}_{k_1^{q+1}, m+(q+1)\nu}, M^{q}_{k_1^{q},m+q\nu})|\leq 4b\Lambda_{m+q\nu}\overset{\mathrm{def}}=\lambda\Lambda_{m+(q-1)\nu}
=\lambda 2^{-\nu(q-1)}\Lambda_m
\end{equation}
and similarly
\begin{equation}\label{20}
|\gamma(M^{q+1}_{k_2^{q+1},m+(q+1)\nu}, M^{q}_{k_2^{q},m+q\nu})|\leq \lambda 2^{-\nu(q-1)}\Lambda_m.
\end{equation}
In the proof of Lemma~\ref{lem2} we obtained the estimates
\begin{equation}\label{21}
|\gamma(M_{k_1^1,m+\nu},M_{k_1^0,m})|\leq a\Lambda_m,
\end{equation}
\begin{equation}\label{22}
|\gamma(M_{k_2^1,m+\nu},M_{k_2^0,m})|\leq a\Lambda_m.
\end{equation}
From \eqref{19}, \eqref{21} and    \eqref{20}, \eqref{22} we find respectively that
\begin{equation}\label{23}
\begin{aligned}
&|\gamma(M^{q+1}_{k_1^{q+1},m+(q+1)\nu}, M_{k_1^0,m}^0)|\leq
(a+\lambda\sum\limits_{j=1}^{q-1}2^{-\nu j})\Lambda_m\\&<
\left(a+\frac18\cdot\frac1{1-2^{-\nu}}\right)
\Lambda_m\leq\left(a+\frac14\right)\Lambda_m
\end{aligned}
\end{equation}
and
\begin{equation}\label{24}
\begin{aligned}
&|\gamma(M^{q+1}_{k_2^{q+1},m+(q+1)\nu}, M_{k_2^0,m}^0)|\leq
\left(a+\frac14\right)\Lambda_m.
\end{aligned}
\end{equation}
From relations \eqref{23}, \eqref{24} and the definition of $a$, we obtain that
 at each next step the points of $L$ added for consideration lie on the arcs
adjacent to the arc $\gamma(M_{k_1,n},M_{k_2,n})$ each of which has length
at most $3b\Lambda_n$. We have
$$\int\limits_{\Omega^0}\cdots=\sum\limits_{q=0}^{\infty}\int\limits_{\Omega^q\setminus \Omega^{q+1}}\cdots \overset{\mathrm{def}}=\sum\limits_{q=0}^\infty I_q.$$
If $M\in\Omega^0\setminus\Omega^1$, then
\begin{equation}\label{25}
\tilde d(M)\geq\underset{M\not\in\Omega^{\ast}_{m+\nu}}\min\dist(M,L)\geq c_8\Lambda_{m+\nu}=c_9\Lambda_m
\end{equation}
with some constant $c_8=c_8(b)$ and $c_9=2^{-\nu}c_8$.
Therefore,
\begin{equation}\label{26}
\begin{aligned}
I_0&\leq c_{10}\cdot\frac1{\Lambda^2_m}\underset{\Omega^0\setminus\Omega^1}\int
\Delta^{\ast}f(N(M),c_6\tilde d(M))dm_3(M)\\ &\leq
c_{10}\cdot\frac1{\Lambda^2_m}\sum\limits_{j=2^{m-n}k_1}^{2^{m-n}k_2}
\underset{\bar B_{\Lambda_m}(M_{jm})}\int\Delta^{\ast}f(N(M),c_6\tilde d(M))dm_3(M).
\end{aligned}
\end{equation}
Since $f\in\Lambda_p^{\alpha_0}(L)$, we have
\begin{equation}\label{27}
\begin{aligned}
\Delta^{\ast}f(N(M),c_6\tilde d(M))&\leq c\cdot(c_6\Lambda_n)^{-\epsilon}
(c_6\tilde d(M))^{\epsilon}\Delta^{\ast}f(M_{k_0n},c_6\Lambda_n)\\ &\leq
c^{\prime}\Lambda^{-\varepsilon}_n\Lambda_m^{\epsilon}\Delta^{\ast}f(M_{k_0n},c_6\Lambda_n).
\end{aligned}
\end{equation}
Then \eqref{26} and \eqref{27} give rise to the estimate
\begin{equation}\label{28}
\begin{aligned}
I_0&\leq c^{\prime\prime}\cdot\frac1{\Lambda_m^2}\sum\limits_{j=2^{m-n}k_1}^{2^{m-n}k_2}
\Lambda_n^{-\epsilon}\Lambda_m^{\epsilon}
\Delta^{\ast}f(M_{k_0n},c_6\Lambda_n)\cdot|B_{\Lambda_m}(\cdot)|\\
&\leq \bar c\Lambda_m2^{m-n}\Lambda_n^{-\epsilon}\Lambda_m^{\epsilon}
\Delta^{\ast}f(M_{k_0n},c_6\Lambda_n)=
c_0^{\prime}\Lambda_n^{1-\epsilon}\Lambda_m^{\epsilon}\Delta^{\ast}f(M_{k_0n},c_6\Lambda_n).
\end{aligned}
\end{equation}
If $M\in \Omega^q\setminus\Omega^{q+1}$, $k\geq1$, then
$\tilde d(M)\leq \Lambda_{m+q\nu}$, $\tilde d(M)\geq d(M)\geq c^{\prime}
\Lambda_{m+(q+1)\nu}$.
Hence
\begin{equation}\label{29}
\begin{aligned}
I_q&\leq c_{10}\frac1{\Lambda_{m+(q+1)\nu}}\sum\limits_{j=k_1^q}^{k_2^q}\;
\underset{\bar B_{m+q\nu}(M_{j,m+q\nu})}\int
\Delta^{\ast} f(N(M), c_6 \tilde d(M))dm_3(M)\\ &\leq
c_{10}\frac1{\Lambda_{m+(q+1)\nu}}\sum\limits_{j=k_1^q}^{k_2^q}\;
\underset{\bar B_{m+q\nu}(M_{j,m+q\nu})}\int
\Delta^{\ast} f(N(M), c_6\Lambda_{m+q\nu})dm_3(M).
\end{aligned}
\end{equation}
As in \eqref{27} we have the estimate
\begin{equation}\label{30}
\Delta^{\ast}f(N(M),c_6\Lambda_{m+q\nu})\leq c^{\prime}\Lambda_n^{-\epsilon}
\Lambda^{\epsilon}_{m+q\nu}\Delta^{\ast}f(M_{k_0n},(c_6+3b)\Lambda_n).
\end{equation}
From \eqref{23} we obtain that
\begin{equation}\label{31}
k_2^q-k_1^q\leq c^{\prime\prime}\Lambda_n\Lambda^{-1}_{m+q\nu}.
\end{equation}
Eqs.\eqref{29} -- \eqref{31}
imply the estimate
\begin{equation}\label{32}
\begin{aligned}
I_q&\leq c_{12}\frac1{\Lambda^2_{m+(q+1)\nu}}\Lambda_n
\Lambda^{-1}_{m+q\nu}\Lambda_n^{-\epsilon}\Lambda_{m+q\nu}^{\epsilon}
\Delta^{\ast}f(M_{k_0n},(c_6+3b)\Lambda_n)\Lambda^3_{m+q\nu}\\
&\leq c_{13}\Lambda_n^{1-\epsilon}\Lambda^{\epsilon}_{m+q\nu}
\Delta^{\ast}f(M_{k_0n},(c_6+3b)\Lambda_n).
\end{aligned}
\end{equation}
Now
\eqref{28} and \eqref{32} imply the estimate
$$\int\limits_{\Omega^0}\cdots\leq   c_{13}\Lambda_n^{1-\epsilon}
\Delta^{\ast}f(M_{k_0n},(c_6+3b)\Lambda_n)\sum\limits_{q=0}^{\infty}
\Lambda^{\epsilon}_{m+q\nu}\leq
c_7\Lambda_n^{1-\epsilon}\Lambda_m^{\epsilon}
\Delta^{\ast}f(M_{k_0n},(c_6+3b)\Lambda_n).$$
\end{proof}
\section{A construction of the approximating function ${v_{\delta}}$ and a representation
of $f_0$ }\label{sec5}
  We construct the function $v_{\delta}$ for $\delta=2^{-n}$;
 for $2^{-n-1}<\delta<2^{-n}$ we put $v_{\delta}=v_{2^{-n}}$.
Define the points $M_{kn}\in L$, $M_{0n}=A$, $M_{2^n,n}=B$
as in the construction of $f_1$.
By a similar reasoning as in \cite {AlexeevaShir-2020}, we choose a $c_{11}\geq 1$ such that
the inequality  $$m_3(\bar B_{c_{11}\Lambda_n}(M_{kn})
\setminus\Omega^{\ast}_{n-2})\geq \frac12m_3(\bar B_{c_{11}\Lambda_n}(M_{kn}))$$
is valid for $0\leq k\leq 2^n$. As is shown in [7], the constant $c_{11}$ depends only on $b$.
      Now we put
$$\beta_{0n}=B_{2\Lambda_n}(M_{0n}),\
\beta_{kn}=B_{2\Lambda_n}(M_{kn})\setminus\bigcup\limits_{\nu=1}^{k-1}
B_{2\Lambda_n}(M_{\nu n}).$$
Geometric considerations imply the inequality
$d(M)\overset{\mathrm{def}}=\dist(M,L)\geq 2^{-n+1}$ for
$M\in B_{c_{11}\Lambda_n}(M_{kn})\setminus \Omega^{\ast}_{n-2}.$
In addition, $d(M)\leq c_{11}2^{-n}|\Lambda|. $
Denoting by $\Delta f_0(M)$ the Laplace operator   for the function $f_0$, we
define the numbers $c_{kn}$ by the equation
\begin{equation}\label{33}
\int\limits_{\beta_{kn}}\Delta f_0(M)dm_3(M)=
c_{kn}\Lambda_n\Delta^{\ast} f(M_{kn},(c_5+2b)\Lambda_n).
\end{equation}
Let $\chi_{kn}$ be the characteristic function of the set
$\bar B_{c_{11}\Lambda_n}(M_{kn})\setminus\Omega^{\ast}_n$
in $\mathbb R^3$. We put
\begin{equation}\label{34}
  \varphi_{kn}(M)=\gamma_{kn}\Lambda^{-2}_n\chi_{kn}(M)
\Delta^{\ast} f(M_{kn},(c_5+2b)\Lambda_n)
\end{equation}
and define the numbers $\gamma_{kn}$ by the equality
\begin{equation}\label{35}
    \int\limits_{\beta_{kn}}\Delta f_0(M)dm_3+
\int\limits_{\mathbb R^3}\varphi_{kn}(M)dm_3=0.
\end{equation}
   The choice of the constant $c_{11}$, Definitions \eqref{33} -- \eqref{35},
and Lemmas \ref{lem1} and \ref{lem2} imply the estimates
$|c_{kn}|\leq c_{12}$, $|\gamma_{kn}|\leq c_{13}$, where $c_{12}$ and $c_{13}$ depend on $f$ and $b$.
 Define
\begin{equation}\label{36}
   \Phi(M)=\sum\limits_{k=0}^{2^n}\varphi_{kn}(M)
\end{equation}
Finally we put
\begin{equation}\label{37}
v_{2^{-n}}(M_0)=-\frac1{4\pi}\int\limits_{\mathbb R^3\setminus
\Omega^{\ast}_n}\frac{\Delta f_0(M)}{||MM_0||}dm_3(M)+
\frac1{4\pi}\int\limits_{\mathbb R^3}\frac{\Phi_n(M)}{||MM_0||}dm_3(M). \end{equation}
  Definitions  \eqref{34}, \eqref{35}, \eqref{36} show that
$v_{2^{-n}}\in H(\Omega_n)$ due to $\Omega_n\subset\Omega^{\ast}_{n-2}.$
  Since  by Proposition~\ref{prop1} we have $f\in \Lambda^{\alpha-\frac1p}(L)$, and the
construction of the continuation of $f_0$ to $\mathbb R^3$  was similar to the
construction of a continuation of a function from the curve $L$ in \cite{AlexeevaShir-2020}, we may apply the reasoning from \cite{AlexeevaShir-2020}
to the special case $\omega(t)=t^{\alpha-\frac1p}(L)$, which gives us the following representation of $f$:
   \begin{equation}\label{38}
            f(M_0)=-\frac1{4\pi}\int\limits_{\mathbb R^3}
\frac{\Delta f_0(M)}{||MM_0||}\,dm_3(M), \ M_0\in L.
   \end{equation}
From \eqref{37} and \eqref{38} we obtain the following expression for the difference
     $v_{2^{-n}}-f$:
\begin{equation}\label{39}
  v_{2^{-n}}(M_0)-f(M_0)=\frac1{4\pi}\int\limits_{\Omega^{\ast}}
\frac{\Delta f_0(M)}{||MM_0||}\,dm_3(M)+\frac1{4\pi}\int\limits_{\mathbb R^3}
\frac{\Phi_0(M)}{||MM_0||}\,dm_3(M).
\end{equation}
\section{Proof of Theorem \ref{thm1}}\label{sec6}
First of all we notice that it is sufficient to
establish the inequalities in Theorem \ref{thm1}
   for the approximating function $v_{2^{-n}}$ with the replacement
of $\max_{2^{-n}}(\cdots)$ by
$\max_{c2^{-n}}(\cdots)$ for some constant $c$.
   Further we will need the following statement.
 \begin{lem}\label{lem3}
  Suppose $f\in\Lambda_p^{\alpha_0}(L)$, the constants
$c>0$ and $\epsilon>0$ are as in \eqref{2},
$c_6\geq 4$, $d(M)=\dist(M,L)$, $M_0\in L$,
$N_0(M)\in L$ is a point such that $||MN_0(M)||=d(M)$.
Then there exists $c_{14}=c_{14}(c,\epsilon,b)>0$ such that
the following estimate is valid:
$$
 \int\limits_{\bar B_{\Lambda_n}(M_0)}
\frac{\Delta^{\ast} f(N_0(M),c_6d(M))}{||MM_0||d^2(M)}\,dm_3(M)\leq
c_{14}\Delta^{\ast} f(M_0,(c_6+5b)\Lambda_n).
$$
    \end{lem}
\begin{proof}
 Consider the inequalities
\begin{equation}\label{41}
\begin{aligned}
\int\limits_{\bar B_{\Lambda_n}(M_0)}
&\frac{\Delta^{\ast} f(N_0(M),c_6d(M))}{||MM_0||d^2(M)}\,dm_3(M)\leq
\sum\limits_{\nu=0}^{\infty}
\frac{2^{\nu+1}}{\Lambda_n}
\underset{2^{-\nu}\bar B_{\Lambda_n}(M_0)
\setminus 2^{-\nu-1}\bar B_{\Lambda_n}(M_0)}\int\cdots\\
&\leq \sum\limits_{\nu=0}^{\infty}\frac{2^{\nu+1}}{\Lambda_n}
\underset{2^{-\nu}\bar B_{\Lambda_n}(M_0)}
\int\frac{\Delta^{\ast}f(N_0(M),c_6d(M))}{d^2(M)}\,dm_3(M).
\end{aligned}
\end{equation}
If $M\in2^{-\nu}\bar B_{\Lambda_n}(M_0)$, $\nu\geq 0$, then
$||M_0N_0(M)||\leq 2b\cdot2^{-\nu}\Lambda_n$.
\end{proof}
Let $s_{\nu}$ be the longest arc of $L$ containing the points
$2^{-\nu}\partial B_{\Lambda_n}(M_0)\cap L
=\partial B_{\Lambda_{n+\nu}}(M_0)\cap L$,
and let $T_{\nu}^{\prime}$, $T_{\nu}^{\prime\prime}$ be  the endpoints of this arc.
Let $\gamma_{\nu}^{\prime}$,  $\gamma_{\nu}^{\prime\prime}\subset L$
be the arcs of length $2b\Lambda_{n+\nu}$
lying outside $s_{\nu}$ and such that $T_{\nu}^{\prime}$ is one of the endpoints of
$\gamma_{\nu}^{\prime}$ and  $T_{\nu}^{\prime\prime}$ is one of the endpoints of
$\gamma_{\nu}^{\prime\prime}$. We put
$S_{\nu}=s_{\nu}\cup \gamma_{\nu}^{\prime}\cup\gamma_{\nu}^{\prime\prime}$
Let $M_{k_1,n+\nu}$ be the point nearest  to $S_{\nu}$
to one side of $S_{\nu}$  and $M_{k_2,n+\nu}$ be the point nearest  to $S_{\nu}$
to the other side of $S_{\nu}$. Then
\begin{equation}\label{42}
|\gamma(M_{k_1,n+\nu}, M_{k_2,n+\nu})|\leq(6b+2)\Lambda_{n+\nu}.
\end{equation}
We choose an $l$ such that $2^{l-4}\geq6b+4$.
Then \eqref{42} implies that every
ball
$2^{- \nu }B_{\Lambda_n}(M_0)$ satisfies the conditions of Lemma \ref{lem1}
with $\Lambda_n$ replaced by $\Lambda_{n+\nu-l}$. Therefore,
\begin{equation}\label{43}
\begin{aligned}
&\underset{2^{-\nu}\bar B_{\Lambda_n}(M_0)}\int
\frac{\Delta^{\ast}f(N_0(M),c_6d(M))}{d^2(M)}dm_3(M)\\
&\leq\underset{\Omega^{\ast}_{n+\nu-l}(2^{\nu-l}k_1^1,2^{\nu-l}k_2^1)}
\int\frac{\Delta^{\ast}f(N_0(M),c_6d(M))}{d^2(M)}dm_3(M)\\
&\leq c_7\Lambda^{1-\epsilon}_{n-l+\nu}\Lambda^{\epsilon}_{n-l+\nu}
\Delta^{\ast}f(M_0,(c_6+5b)\Lambda_{n-l})\\&=c_7\Lambda_{n-l+\nu}\Delta^{\ast}
f(M_0,(c_6+5b)\Lambda_{n-l+\nu})\\
&\leq c_{15}\Lambda_{n+\nu}
\left(\frac{\Lambda_{n+\nu}}{\Lambda_n}\right)^{\epsilon}
\Delta^{\ast}f(M_0,(c_6+5b)\Lambda_n).
\end{aligned}
\end{equation}
From \eqref{41} and \eqref{43} we find that
$$
\begin{aligned}
&\int\limits_{\bar B_{\Lambda_n}(M_0)}
\frac{\Delta^{\ast}f(N_0(M),c_6d(M))}{||MM_0||d^2(M)}\,dm_3(M)\\
&\leq c_{15}\sum\limits_{\nu=0}^{\infty}
\frac{2^{\nu+1}}{\Lambda_n}\Lambda_{n+\nu}
\left(\frac{\Lambda_{n+\nu}}{\Lambda_n}\right)^{\epsilon}
\Delta^{\ast}f(M_0,(c_6+3b)\Lambda_n)\\
&\leq c_{15}^{\prime}\sum\limits_{\nu=0}^{\infty}
2^{-\nu\epsilon}\Delta^{\ast}f(M_0,(c_6+5b)\Lambda_n)=
c_{14}\Delta^{\ast}f(M_0,(c_6+5b)\Lambda_n).
\end{aligned}
$$
The lemma is proved.
Now we choose an arbitrary measurable function $K(M)$, $M\in L$,
such that $||MK(M)\leq\Lambda_n$, $K(M)\in L$, and define the  functions
$G_m(M)$, where $m$ is an integer such that $|m|\leq 2^n$ as follows.
For $M\in\gamma(M_{kn},M_{k+1,n})$, where $0\leq k\leq 2^n-1$, denote by
$Q_{kn}(M)$ the set of all indices $l$, $0\leq l\leq 2^n$, for which
the relation
\begin{equation}\label{45}
\beta_{ln}\cap\bar B_{3\Lambda_n}(M_{kn})\neq\emptyset
\end{equation}
and put
\begin{equation}\label{46}
G_0(M)\overset{\mathrm{def}}=
\sum\limits_{l\in Q_{kn}(M)}\left(\frac1{4\pi}\int\limits_{\beta_{ln}}
\frac{\Delta f_0(P)}{||PM||}\,dm_3(P)+
\frac1{4\pi}\int\limits_{\mathbb R^3}\frac{\varphi_{ln}(P)}{||PM||}\,dm_3(P)\right).
\end{equation}
For $m\neq0$, $|m|\leq2^n$, the functions $G_m(M)$ are defined as follows.
We put
\begin{equation}\label{46}
G_m(M)=\begin{cases}0 \  \text{if} \   m+k\in Q_{kn}(M),\ \text{or}\
\beta_{m+k,n}=\emptyset,\ \text{or}\ m+k\not\in[0,2^n],\\
\displaystyle\frac1{4\pi}\int\limits_{\beta_{m+k,n}}
\frac{\Delta f_0(P)}{||PM||}dm_3(P)+\frac1{4\pi}
\int\limits_{\mathbb R^3}\frac{\varphi_{m+k,n}(P)}{||PM||}dm_3(M) \
\text{otherwise}.
\end{cases}
\end{equation}
From \eqref{39}, \eqref{45}, \eqref{46} we get the equality
\begin{equation}\label{47}
v_{2^{-n}}(M)-f(M)=\sum\limits_{m=-2^n}^{2^n}G_m(M),\ M\in L.
\end{equation}
Now \eqref{47} implies
\begin{equation}\label{48}
\begin{aligned}
\left(\int\limits_L|v_{2^{-n}}(K(M))-f(K(M))|^p\,dm_1(M)\right)^{\frac1p}\\
\leq\sum\limits_{m=-2^n}^{2^n}
\left(\int\limits_L|G_m(K(M))|^p\,dm_1(M)\right)^{\frac1p}\overset{\mathrm{def}}=
\sum\limits_{m=-2^n}^{2^n}I_m.
\end{aligned}
\end{equation}
Let us estimate the term $I_0$.
Lemma~\ref{lem1} yields the estimate
\begin{equation}\label{49}
\begin{aligned}
I_0^p&=\sum\limits_{k=0}^{2^n-1}\
\int\limits_{\gamma(M_{kn},M_{k+1,n})}G_0^p(K(M))\,dm_1(M)\\
&\leq c\sum\limits_{k=0}^{2^n-1}\ \int\limits_{\gamma(M_{kn},M_{k+1,n})}
\left(\sum\limits_{l\in Q_{kn}(M)}\ \int\limits_{\beta_{ln}(M)}
\frac{\Delta f_0(P)}{||PK(M)||}\,dm_3(P)\right)^pdm_1(M)\\
&+c\sum\limits_{k=0}^{2^n-1}\ \int\limits_{\gamma(M_{kn},M_{k+1,n})}
\left(\sum\limits_{l\in Q_{kn}(M)}\ \int\limits_{\mathbb R^3}
\frac{\varphi_{ln}(P)}{||PK(M)||}\,dm_3(P)\right)^pdm_1(M)\\
&\leq c\sum\limits_{k=0}^{2^n-1}\ \int\limits_{\gamma(M_{kn},M_{k+1,n})}
\left(\int\limits_{4B_{\Lambda_n}(M_{kn})}
\frac{\Delta^{\ast}f(N(P),c_5d(P))}{||PK(M)||d^2(P)}dm_3(P)\right)^pdm_1(M)\\
&+c\sum\limits_{k=0}^{2^n-1}\ \int\limits_{\gamma(M_{kn},M_{k+1,n})}
\left(\sum\limits_{l\in Q_{kn}(M)}\ \int\limits_{\mathbb R^3}
\frac{\varphi_{ln}(P)}{||PK(M)||}\,dm_3(P)\right)^pdm_1(M),
\end{aligned}
\end{equation}
where $N(P)$ in \eqref{49} is a point on $L$
for which $||PN(P)||=d(P)$, $d(P)=\dist(P,L)$.
Replacing $c_5+3b$, $c_5+5b$ by $c_{15}$, $c_{16}$, etc.,
we obtain from Lemma~\ref{lem3} the inequality
\begin{equation}\label{50}
\begin{aligned}
\int\limits_{4\bar B_{\Lambda_n}(M_{kn})}
\frac{\Delta^{\ast}f(N(P),c_5d(P))}{||PK(M)||d^2(P)}dm_3(P)
\leq c\Delta^{\ast}f(M_{kn},c_{15}\Lambda_n).
\end{aligned}
\end{equation}
From \eqref{50}
we get the inequality
\begin{equation}\label{51}
\begin{aligned}
&\sum\limits_{k=0}^{2^n-1}\ \int\limits_{\gamma(M_{kn},M_{k+1,n})}
\left(\int\limits_{4B_{\Lambda_n}(M_{kn})}
\frac{\Delta^{\ast}f(N(P),c_5d(P))}{||PK(M)||d^2(P)}dm_3(P)\right)^pdm_1(P)\\
&\leq c\sum\limits_{k=0}^{2^{n-1}}\int\limits_{\gamma(M_{kn},M_{k+1,n})}
(\Delta^{\ast}f(M_{kn},c_{15}\Lambda_n))^pdm_1(P)\\
&\leq c\int\limits_L(\Delta^{\ast}f(M,c_{16}\Lambda_n))^pdm_1(P)\leq c2^{-2\alpha p}.
\end{aligned}
\end{equation}
From  \eqref{33} and the inequality
$|\gamma_{kn}|\leq c_{13}$, we obtain an estimate similar to \eqref{51} for the second summand in \eqref{49}, and as a result we get the inequality
\begin{equation}\label{51}
I_0\leq c2^{-n\alpha}.
\end{equation}
Estimating $I_m$ for $m\neq0$, we take into account that if
$k+m\not\in Q_{kn}(M)$  and $\beta_{k+m,n}\neq\emptyset$, then
for $P\in\beta_{k+m,n}$ we have
\begin{equation}\label{52}
\begin{aligned}
&||PM||\geq3\Lambda_n, \ ||PM_{k+m,n}||\leq2\Lambda_n\leq\frac23||PM||,\\
&||MM_{k+m,n}||\leq||PM||+||PM_{k+m,n}||\leq\frac53||PM||.
\end{aligned}
\end{equation}
At the same time due to the $b$-chord-arc condition on  $L$ and the definition of the set
$Q_{kn}(M)$ we have
\begin{equation}\label{53}
||MM_{k+m,n}||\geq\frac1b|\gamma(M,M_{k+m,n})|\geq\frac1b(|m|-1)\Lambda_n.
\end{equation}
Taking into account that $|m|\geq2$ whenever $k+m\not\in Q_{kn}^0$,
we obtain from \eqref{52} and \eqref{53} that
$||PM||\geq(3/10 b)\Lambda_n|m|$. From \eqref{34} and \eqref{45}
we obtain
\begin{equation}\label{54}
\begin{aligned}
G_m(M)&=\frac1{4\pi}\int\limits_{\beta_{k+m,n}}
\frac{\Delta f_0(P)}{||PM||}dm_3(P)+\frac1{4\pi}\int_{\mathbb R^3}
\frac{\varphi_{k+m,n}(P)}{||PM||}dm_3(P)\\&-
\frac1{4\pi}\int\limits_{\beta_{k+m,n}}\frac{\Delta f_0(P)}{||M_{k+m,n}M||}
dm_3(P)-\frac1{4\pi}\int\limits_{\mathbb R^3}\frac{\varphi_{k+m,n}(P)}{||M_{k+m,n}M||}dm_3(P)\\
&=\frac1{4\pi}\int\limits_{\beta_{k+m,n}}\Delta f_0(P)\left(\frac1{||PM||}-\frac1{||M_{k+m,n}M||}\right)dm_3(P)\\&+
\frac1{4\pi}\int\limits_{\mathbb R^3}\varphi_{k+m,n}(P)
\left(\frac1{||PM||}-\frac1{||M_{k+m,n}M||}\right)dm_3(P).
\end{aligned}
\end{equation}
Due to $||PM||\geq(3/10 b)\Lambda_n|m|$ and $||M_{k+m,n}M||\geq(1/b)|m|\Lambda_n$
we have
$$\left|\frac1{||PM||}-\frac1{||M_{k+m,n}M||}\right|\leq c\frac 1{m^2\Lambda_n}.$$
Hence by \eqref{16} and Lemma~\ref{2} we obtain the estimate
\begin{equation}\label{55}
\begin{aligned}
\left|\frac1{4\pi}\int\limits_{\beta_{k+m,n}}\Delta f_0(P)
\left(\frac1{||PM||}-\frac1{||M_{k+m,n}M||}\right)dm_3(P)\right|\\
\leq c\frac1{m^2\Lambda_n}\
\int\limits_{\beta_{k+m,n}}
\frac{\Delta^{\ast}f(N(P),c_5d(P))}{d^2(P)}dm_3(P)\\
\leq  c\frac1{m^2\Lambda_n}\
\int\limits_{B_{2\Lambda_{n-2}}(M_{k+m,n})}
\frac{\Delta^{\ast}f(N(P),c_5d(P))}{d^2(P)}dm_3(P)\\
\leq c\frac1{m^2\Lambda_n}\Lambda_n\Delta^{\ast}f(M_{k+m,n},c_{16}\Lambda_n)
\leq c\frac1{m^2}\Delta^{\ast}f(M_{k+m,n},c_{16}\Lambda_n).
\end{aligned}
\end{equation}
Similarly we obtain the inequality
\begin{equation}\label{56}
\begin{aligned}
\frac1{4\pi}\int\limits_{\mathbb R^3}\varphi_{k+m,n}(P)
\left(\frac1{||PM||}-\frac1{||M_{k+m,n}M||}\right)dm_3(P)
\leq  c\frac1{m^2}\Delta^{\ast}f(M_{k+m,n},c_{16}\Lambda_n).
\end{aligned}
\end{equation}
As a result, we have the following estimate from \eqref{54} -- \eqref{56}:
\begin{equation}\label{57}
G_m(M)\leq c\frac1{m^2}\Delta^{\ast}f(M_{k+m,n},c_{16}\Lambda_n).
\end{equation}
In the case where $k+m\not\in [0,2^n]$, both  sides of \eqref{57}
are zero and inequality~\eqref{57} is valid for all $M\in\gamma(M_{kn},M_{k+1,n})$.
Since the function $K(M)$ satisfies the condition
$||MK(M)||\leq\Lambda_n$, estimate \eqref{57}
implies the inequality
\begin{equation}\label{58}
G_m(K(M))\leq c\frac1{m^2}\Delta^{\ast}f(M_{k+m,n},c_{17}\Lambda_n)
\end{equation}
with the constant $c_{17}=c_{16}+2b$.
Now for $|m|\geq2$ it follows from \eqref{58}
that
$$
\begin{aligned}
I_m^p=\int\limits_{L}|G_m(K(M))|^p\,dm_1(M)&\leq
\sum\limits_{k=0}^{2^n-1}c \int\limits_{\gamma(M_{kn},M_{k+1,n})}
\frac{1}{m^{2p}}(\Delta^{\ast}f(M_{k+m,n},c_{17}\Lambda_n))^p\,dm_1(P)\\
&\leq c\frac1{m^{2p}}\int\limits_L(\Delta^{\ast}f(M,c_{18}\Lambda_n))^p\,dm_1(P)
\leq c\frac1{m^{2p}}
2^{-n\alpha p},
\end{aligned}
$$
hence
\begin{equation}\label{60}
I_m\leq c\frac1{m^2}2^{-n\alpha}.
\end{equation}
From estimates \eqref{48}, \eqref{51}, and \eqref{60}
we find that
\begin{equation}\label{61}
\left(\int\limits_L|v_{2^{-n}}(K(M))-f(K(M))|^p\,dm_1(M)\right)^{\frac1p}\\
\leq  c\sum\limits_{m=-2^n}^{2^n}\frac{1}{|m|+1)^2}2^{-n\alpha}\leq c2^{-n\alpha}.
\end{equation}
Since   $K(M)$ is an arbitrary function measurable with respect to $m_1$-measure on $L$
with the condition $||MK(M)||\leq \Lambda_n$, we see that \eqref{61} implies \eqref{3}.
Now we proceed to the proof of relation \eqref{4}.
\begin{lem}\label{lem4}
Let
\begin{equation}\label{62}
U_n(M)=\frac1{4\pi}\int\limits_{\mathbb R^3}\frac{\Phi_n(P)}{||PM||}\,dm_3(P).
\end{equation}
Then
\begin{equation}\label{63}
\left(\int\limits_{L}(\grad^{\ast}_{2^{-n}}U_n(M))^p\,dm_1(M)\right)^{\frac1p}\leq c2^{n(1-\alpha)}.
\end{equation}
\end{lem}
\begin{proof}
Equations \eqref{33} and \eqref{35}
imply the equality
\begin{equation}\label{64}
\begin{aligned}
U_n(M)&=\frac1{4\pi}\int\limits_{\mathbb R^3\setminus \Omega^{\ast}_n}
\frac{\sum_{k=0}^{2^n}\varphi_{kn}(P)}{||PM||}dm_3(P)\\
&=\frac1{4\pi}\sum\limits_{k=0}^{2^n}
\gamma_{kn}\Lambda_n^{-2}
\int\limits_{\bar B_{c_{11}\Lambda_n}(M_{kn})\setminus\Omega^{\ast}_{n-2}}
\frac{\Delta^{\ast}f(M_{kn}(c_5+2b)\Lambda_{n-2})}{||PM||}dm_3(P).
\end{aligned}
\end{equation}
If $M\in \gamma(M_{jn},M_{j+1,n})$, $0\leq j\leq 2^n$, then
for all $P\in \bar B_{c_{11}\Lambda_n}(M_{kn})\setminus\Omega^{\ast}_n$
we have
$||PM||\geq c(|k-j|+1)\Lambda_n$. Therefore \eqref{64} and the property
$|\gamma_{kn}|\leq c$ yield the following estimate:
\begin{equation}\label{65}
\begin{aligned}
\grad^{\ast}_{2^{-n}}U_n(M)&\leq
c\sum\limits_{k=0}^{2^n}\Lambda_n^{-2}\frac1{(|k-j|+1)^2\Lambda_n^2}
\Lambda_n^3
\Delta^{\ast}f(M_{kn},(c_5+2b)\Lambda_n)\\
&=c\sum\limits_{k=0}^{2^n}\frac1{(|k-j|+1)^2\Lambda_n}
\Delta^{\ast}f(M_{kn},(c_5+2b)\Lambda_n).
\end{aligned}
\end{equation}
Applying the reasoning used for estimating the function $G_m$
to the last sum in \eqref{65}, we obtain that
$$\begin{aligned}
&\int\limits_L(\grad^{\ast}_{2^{-n}}U_n(M))^pdm_1(M))^{\frac1p}\\
\leq c\sum\limits_{m=1}^{2^n+1}\frac1{m^2\Lambda_n}
&\left(\int\limits_L(\Delta^{\ast}f(M,c_{18}\Lambda_n))^pdm_1(M)\right)
^{\frac1p}\leq c2^{n(1-\alpha)},
\end{aligned}
$$
which proves \eqref{63}.
The lemma is proved.
\end{proof}
Let
$$V_n(M)=-\frac1{4\pi}\int\limits_{\mathbb R^3\setminus \Omega^{\ast}_n}
\frac{\Delta f_0(P)}{||PM||}dm_3(P), \quad M\in L.$$
\begin{lem}\label{lem5}
For the function $V_n$ harmonic in $\Omega_{2^{-n}}$ the following
inequality is valid:
\begin{equation}\label{66}
\int\limits_L(\grad^{\ast}_{2^{-n}}V_n(M))^pdm_1(M)\leq c2^{np(1-\alpha)}.
\end{equation}
\end{lem}
\begin{proof}
Notice that $\supp f_0\subset \Omega^{\ast}_0=\bar B_{2|\Lambda|}(A)$ and define
\begin{equation}\label{67}
W_{n\nu}(M)=-\frac1{4\pi}\int\limits_{\Omega^{\ast}_{\nu-1}\setminus
\Omega_{\nu}^{\ast}}\frac{\Delta f_0(P)}{||PM||}dm_3(P),\quad
1\leq\nu\leq n.
\end{equation}
Then it follows from \eqref{66} and \eqref{67} that
$V_n(M)=\sum\limits_{\nu=1}^{n}W_{n\nu}(M).$
Therefore,
\begin{equation}\label{68}
\begin{aligned}
\left(\int\limits_L(\grad^{\ast}_{2^{-n}}V_n(M))^pdm_1(M)\right)^{\frac1p}
\leq \left(\int\limits_L(\sum\limits_{\nu=1}^n
\grad^{\ast}_{2^{-n}}W_{n\nu}(M))^pdm_1(M)\right)^{\frac1p}\\
\leq \sum\limits_{\nu=1}^n\left(\int\limits_L(
\grad^{\ast}_{2^{-n}}W_{n\nu}(M))^pdm_1(M)\right)^{\frac1p}.
\end{aligned}
\end{equation}
Let
$M\in \gamma(M_{k,\nu-1},M_{k+1,\nu-1}),\ 2\leq\nu\leq n,\ 0\leq k\leq 2^{\nu-1}-1.$
Proceeding in the same way as in the final step of the proof of \eqref{3}, we
denote by $Q_{k,\nu-1}(M)$  the set of indices $l$, $0\leq l\leq 2^{\nu-1}$, for which $\beta_{l,\nu-1}\cap\bar B_{3\Lambda_{\nu-1}}(M_{k,\nu-1})\neq\emptyset$
and define
\begin{equation}\label{69}
G_{0\nu}(M)=\sum\limits_{l\in Q_{k,\nu-1}(M)}
\left(-\frac1{4\pi}
\int\limits_{\beta_{l,\nu-1}\setminus\Omega^{\ast}_{\nu}}
\frac{\Delta f_0(P)}{||PM||}dm_3(P)\right).
\end{equation}
For $m\neq0$, $|m|\leq2^{\nu-1}$,
$k+m\not\in Q_{k,\nu-1}(M)$, and $\beta_{k+m,\nu-1}\setminus
\Omega^{\ast}_{\nu}\neq\emptyset$,
the functions $G_{m\nu}(M)$ are defined by the equation
\begin{equation}\label{70}
G_{m\nu}(M)=-\frac1{4\pi}\int\limits_{\beta_{k+m,\nu-1}\setminus\Omega^{\ast}_{\nu}}
\frac{\Delta f_0(P)}{||MP||}\,dm_3(P).
\end{equation}
If the above conditions are not satisfied or $k+m\not\in[0,2^{\nu-1}]$,
then we put $G_{m\nu}(M)=0$.
From \eqref{67}, \eqref{68}, and \eqref{70} we obtain that
\begin{equation}\label{71}
W_{n\nu}(M)=\sum\limits_{m=-2^{\nu-1}}^{2^{\nu-1}}G_{m\nu}(M),
\end{equation}
and it follows from \eqref{71} that
$$\grad^{\ast}_{2^{-n}}W_{n\nu}(M)\leq
\sum\limits_{m=-2^{\nu-1}}^{2^{\nu-1}}\grad^{\ast}_{2^{-n}}G_{m\nu}(M),$$
which implies the inequality
\begin{equation}\label{72}
\left(\int\limits_L(\grad^{\ast}_{2^{-n}}
W_{n\nu}(M))^pdm_1(M)\right)^{\frac1p}
\leq
\sum\limits_{m=-2^{\nu-1}}^{2^{\nu-1}}\left(\int\limits_L(\grad^{\ast}_{2^{-n}}
G_{m\nu}(M))^pdm_1(M)\right)^{\frac1p}.
\end{equation}
If $P\in(\beta_{l,\nu-1}\cap\bar B_{3\Lambda_{\nu-1}}(M_{k,\nu-1}))\setminus \Omega^{\ast}_{\nu}$,\
$M\in \gamma(M_{k,\nu-1}, M_{k+1,\nu-1})$, then $||MP||\geq c\Lambda_{\nu}$,
$d(P)\geq \Lambda_{\nu}$,\ $d(P)\leq 2\Lambda_{\nu-1}=4\Lambda_{\nu}$.
It follows from Lemma~\ref{lem1} and \eqref{69} that
$$\grad^{\ast}_{2^{-n}}G_{0\nu}(M)\leq
c\frac1{\Lambda_{\nu}^2}\frac{\Delta^{\ast}f(M_{k,\nu-1},8c_5\Lambda_{\nu})}
{\Lambda_{\nu}^2}\Lambda_{\nu}^3=
c\frac{\Delta^{\ast}f(M_{k,\nu-1},8c_5\Lambda_{\nu})}{\Lambda_{\nu}}.
$$
Therefore,
\begin{equation}\label{73}
\begin{aligned}
\left(\int\limits_L(\grad^{\ast}_{2^{-n}}
G_{0\nu}(M))^pdm_1(M)\right)^{\frac1p}\\
\leq\frac c{\Lambda_{\nu}}\left(\int\limits_L(
\Delta^{\ast}f(M,10c_5\Lambda_{\nu}))^pdm_1(M)\right)^{\frac1p}
\leq c2^{\nu(1-\alpha)}.
\end{aligned}
\end{equation}
If $G_m(M)\neq0$, $m\neq0$, then for $P\in \beta_{k+m,\nu-1}\setminus\Omega^{\ast}_{\nu}$  we have
$d(P)\geq\Lambda_{\nu}$, $|MP|\geq c|m|\Lambda_{\nu}$,
$d(P)\leq4\Lambda_{\nu}$. Again by \eqref{70} and Lemma~\ref{lem1}
we get
\begin{equation}\label{74}
\begin{aligned}
\grad^{\ast}_{2^{-\nu}}G_{m\nu}(M)\leq c\frac1{m^2\Lambda_{\nu}^2}
\frac{\Delta^{\ast}f(M_{k+m,\nu-1},8c_5\Lambda_{\nu})}{\Lambda_{\nu}^2}
\Lambda_{\nu}^3=c\frac{\Delta^{\ast}f(M_{k+m,\nu-1},8c_5\Lambda_{\nu})}
{m^2\Lambda_{\nu}.}
\end{aligned}
\end{equation}
From \eqref{74} we obtain the estimate
\begin{equation}\label{75}
\begin{aligned}
&\left(\int\limits_L(\grad^{\ast}_{2^{-n}}
G_{m\nu}(M))^pdm_1(M)\right)^{\frac1p}\\
\leq c \frac1{m^2\Lambda_{\nu}}&\left(\int\limits_L(\Delta^{\ast}
f(M,10c_5\Lambda_{\nu}))^p
dm_1(M)\right)^{\frac1p}
\leq c\frac1{m^2}2^{\nu(1-\alpha)}.
\end{aligned}
\end{equation}
From \eqref{72}, \eqref{73}, and \eqref{75} we see that
\begin{equation}\label{76}
\begin{aligned}
 \left(\int\limits_{L}(\grad^{\ast}_{2^{-n}}
W_{n\nu}(M))^pdm_1(M)\right)^{\frac1p}
\leq c\left(2^{\nu(1-\alpha)}+\sum\limits_{m=1}^{2^{\nu-1}}
\frac1{m^2}2^{\nu(1-\alpha)}\right)\leq c2^{\nu(1-\alpha)}.
\end{aligned}
\end{equation}
Now \eqref{68} and \eqref{76} imply the inequality
$$
 \left(\int\limits_{L}(\grad^{\ast}_{2^{-n}}
V_{n}(M))^pdm_1(M)\right)^{\frac1p} \leq
c\sum\limits_{\nu=1}^{n}2^{\nu(1-\alpha)}\leq c2^{n(1-\alpha)},
$$
which is equivalent to \eqref{66}.
Lemma~\ref{lem5} is proved.
Inequality~\ref{4} follows from the equality
$v_{2^{-n}}(M)=V_n(M)+U_n(M)$ and Lemmas~\ref{4} and \ref{5}.
Theorem~\ref{thm1} is proved.
\section{Conclusion}\label{sec7}
It is of interest to find out whether Theorem~1 is true for the entire class
$\Lambda_p^{\alpha}$ under the condition $0<\alpha<1$, $p>1/\alpha$.
\end{proof}



\begin{thebibliography}{1}
\expandafter\ifx\csname url\endcsname\relax
  \def\url#1{\texttt{#1}}\fi
\expandafter\ifx\csname urlprefix\endcsname\relax\def\urlprefix{URL }\fi
\expandafter\ifx\csname href\endcsname\relax
  \def\href#1#2{#2} \def\path#1{#1}\fi

\bibitem{Timan-Book-1963}
A.~F. Timan, Theory of Approximation of Functions of a Real Variable, 1st
  Edition, Vol.~34, Pergamon, 1963.

\bibitem{Motornyi-1971}
V.~P. Motornyi, Approximation of functions by algebraic polynomials in the
  ${L}^p$ metric, Math. USSR-Izv 5~(4) (1971) 889–--914.

\bibitem{Potapov-1977}
M.~K. Potapov, The structural characteristics of the classes of functions with
  a given order of best approximation, Trudy Mat. Inst. Steklov. 134 (1977)
  295–314.

\bibitem{NevaiYu-1994}
P.~Nevai, X.~{Yuan}, Mean convergence of {H}ermite interpolation, Journal of
  Approximation Theory 77~(3) (1994) 282--304.

\bibitem{Dynkin-1983}
E.~M. Dyn'kin, A constructive characterization of classes of
  {S}.~{L}.~{S}obolev and {O}.~{V}.~{B}esov, Proc. Steklov Inst. Math. 155
  (1983) 39--–74.

\bibitem{Andrievskii95}
V.~V. Andrievskii, V.~V. {Maimeskul}, Constructive description of certain
  classes of functions on quasismooth arcs, Russian Academy of Sciences.
  Izvestiya Mathematics 44~(1) (1995) 193--206.
\newblock \href
  {https://doi.org/http://dx.doi.org/10.1070/IM1995v044n01ABEH001589}
  {\path{doi:http://dx.doi.org/10.1070/IM1995v044n01ABEH001589}}.

\bibitem{AlexeevaShir-2020}
T.~A. Alexeeva, N.~A. {Shirokov}, Constructive description of {H}\"{o}lder-like
  classes on an arc in ${ R}^3$ by means of harmonic functions, Journal of
  Approximation Theory 249 (2020) 105308.

\bibitem{Andrievskii88}
V.~V. Andrievskii, On approximation of functions by harmonic polynomials,
  Mathematics of the USSR-Izvestiya 30~(1) (1988) 1--13.
\newblock \href
  {https://doi.org/http://dx.doi.org/10.1070/IM1988v030n01ABEH000989}
  {\path{doi:http://dx.doi.org/10.1070/IM1988v030n01ABEH000989}}.

\bibitem{Stein-B-1970}
E.~M. Stein, Singular Integrals and Differentiability Properties of Functions,
  Princeton Univ. Press, 1970.

\end{thebibliography}
\bibliographystyle{elsarticle-num}

\end{document}